\documentclass{amsart}
\usepackage{euscript}           % Nice script font.
\usepackage{supertabular}
\usepackage{amsmath,amsthm}     % Handy math stuff, theorem environments.
\usepackage{amssymb}            % Fancy math symbols.

\newcommand{\Q}{{\mathbb Q}}

\DeclareMathOperator*{\colim}{colim}

\newcommand {\Z}{{\mathbb Z}}

\newcommand{\lra}{\longrightarrow}              % long right arrow
\newfont{\german}       {eufm10 at 12pt}
\DeclareMathOperator{\Hom}{Hom} \DeclareMathOperator{\Ext}{Ext}

\numberwithin{equation}{section}
\newtheorem{thm}[equation]{Theorem}
\newcounter{numerierer}
\newcounter{leer}

\newtheorem{defn}[equation]{Definition}
\newtheorem{prop}[equation]{Proposition}

\newtheorem{cor}[equation]{Corollary}
\newtheorem{lemma}[equation]{Lemma}

\theoremstyle{definition}  % Bold headings and Roman body text.

\newtheorem{note}[equation]{Note}
\newtheorem{remark}[equation]{Remark}
\usepackage[frame,matrix,curve, arrow]{xy}
\xymatrixcolsep{1.9pc}                          % Adjust size of diagrams.
\xymatrixrowsep{1.9pc}
\newdir{ >}{{}*!/-5pt/\dir{>}}                  % Make better tailed arrows.
\newdir{ |}{{}*!/-5pt/\dir{|}}                  % Make better mapsto arrows.

\newcommand{\mc}[1]{\mathcal{#1}}
\newcommand{\mb}[1]{\mathbb{#1}}

\newcommand{\ul}[1]{\underline{#1}}
\newcommand{\ZZ}{\mathbb{Z}}

\newcommand{\QQ}{\mathbb{Q}}

\newcommand{\FF}{\mathbb{F}}

\newcommand{\MS}{\mathbb{S}}

\DeclareMathOperator*{\Tot}{Tot}
\DeclareMathOperator{\TMF}{TMF}
\DeclareMathOperator{\Gal}{Gal}
\DeclareMathOperator{\ord}{ord}

\begin{document}

\title{$\beta$-family congruences
  and the $f$-invariant}
\author[Mark Behrens and Gerd Laures]{Mark Behrens$\sp 1$ and Gerd Laures}

\address{
Department of Mathematics, MIT, 77 Massachusetts Avenue
Cambridge, Ma 02139-4307, USA \newline\indent
Fakult\"at f\"ur Mathematik,  Ruhr-Universit\"at Bochum, NA1/66,
  D-44780 Bochum, Germany}
%\email{gerd@laures.de }
\date{\today}
\begin{abstract}
In previous work, the authors have each introduced methods for studying the
$2$-line of the $p$-local Adams-Novikov spectral sequence in terms of the
arithmetic of modular forms. 
We give the precise relationship between the congruences of modular forms
introduced by the first author with the $Q$-spectrum
and the $f$-invariant of the second author. This relationship 
enables us to refine
the target group of the $f$-invariant in a way which makes it more manageable 
for computations.
\end{abstract}
\maketitle

\footnotetext[1]{The first author was supported by the NSF grant DMS-0605100, 
the Sloan Foundation, and DARPA.}

\section{Introduction}\label{sec:intro}

In \cite{MR0198470}, J.F.~Adams studied the image of the $J$-homomorphism
$$ J: \pi_{t}(SO) \rightarrow \pi^S_{t} $$
by introducing a pair of invariants
\begin{gather*}
d = d_t: \pi_t^S \rightarrow \pi_t K, \\ 
e = e_t: \ker(d_t) \rightarrow \Ext^{1,t+1}_{\mc{A}}(K_*, K_*) 
\end{gather*}
where $\mc{A}$ is a certain abelian category of graded abelian groups with
Adams operations.  (Adams also studied analogs of $d$ and $e$ using real
$K$-theory, to more
fully detect
$2$-primary phenomena.)
In order to facilitate the study of the $e$-invariant,
Adams used the Chern character to provide a monomorphism
$$ \theta_S: \Ext^{1,t+1}_{\mc{A}}(K_*, K_*) \hookrightarrow \QQ/\ZZ. $$
Thus, the $e$-invariant may be regarded as taking values in $\QQ/\ZZ$.
Furthermore, he showed that for $t$ odd, and $k = (t+1)/2$, 
the image of $\theta_S$ is the cyclic group of
order $\mathrm{denom} (B_{k}/2k)$, where $B_k$ is the $k$th Bernoulli
number.

The $d$ and $e$-invariants detect the $0$ and $1$-lines of the
Adams-Novikov spectral sequence (ANSS).  In \cite{MR1660325}, the second
author studied an invariant
$$ f: \ker(e_t) \rightarrow \Ext^{2,
t+2}_{\TMF_*\TMF[\frac{1}{6}]}(\TMF[\tfrac{1}{6}]_*,
\TMF[\tfrac{1}{6}]_*) $$
which detects the $2$-line of the ANSS for $\pi_*^S$ away from the primes
$2$ and $3$.  He furthermore used H.~Miller's elliptic character to show
that, if $t$ is even and $k = (t+2)/2$, there is a monomorphism
$$ \iota^2: \Ext^{2, t+2}_{\TMF_*\TMF[\frac{1}{6}]}(\TMF[\tfrac{1}{6}]_*,
\TMF[\tfrac{1}{6}]_*) \hookrightarrow 
D_\QQ/(D_{\ZZ[\frac{1}{6}]} + (M_0)_\QQ + (M_k)_\QQ),
$$
where $D$ is Katz's ring of divided congruences and $M_k$ is the space of
weight $k$ modular forms of level $1$ meromorphic at the cusp.  It is
natural to ask for a description of the image of the
map $\iota^2$ in arithmetic terms.

\begin{remark}
In \cite{MR1660325}, the second author works with more general congruence
subgroups $\Gamma  \subseteq SL_2(\ZZ)$ and associated cohomology theories
$E^\Gamma$ which also lead to results for the primes 2 and 3.  
The spectrum $\TMF$ is just the spectrum
$E^{SL_2(\ZZ)}$ when 6 is inverted.  
In this paper we shall not be considering the $f$
invariant associated to more general
congruence subgroups $\Gamma$ and 6 shall always be a unit.
\end{remark}

Attempting to generalize the $J$ fiber-sequence
$$ J \rightarrow KO_p \xrightarrow{\psi^\ell-1} KO_p $$
the first author introduced a ring spectrum $Q(\ell)$ built from a length two
$\TMF_p$-resolution.
In \cite[Thm.~12.1]{BM08}, it was shown that for $p \ge 5$, 
the elements $\beta_{i/j,k} \in
(\pi_*^S)_{p}$ of \cite{MillerRavenelWilson} are detected in the Hurewicz
image of $Q(\ell)$.  This gives rise to the association of a modular form
$f_{i/j,k}$ to each element $\beta_{i/j,k}$.  Furthermore, the forms
$f_{i/j,k}$ are characterized by certain
arithmetic conditions.

The purpose of this paper is to prove that the $f$-invariant of
$\beta_{i/j,k}$ is given by the formula
$$ f(\beta_{i/j,k}) = \frac{f_{i/j,k}}{p^k E_{p-1}^j} \qquad
\text{(Theorem~\ref{thm:mainthm})}. $$
In particular, since the $2$-line of the ANSS is generated by the elements
$\beta_{i/j,k}$, the $p$-component of the 
image of the map $\iota^2$ is characterized by the arithmetic conditions
satisfied by the elements $f_{i/j,k}$.  

J.~Hornbostel and N.~Naumann \cite{HornbostelNaumann} computed the $f$ invariant of the
elements $\beta_{i/1,1}$ in terms of Katz's Artin-Schreier generators of the ring of 
$p$-adic modular forms.  
While their result is best suited to describe $f$-invariants of infinite families, it is
difficult to explicitly get one's hands on their output.
Direct computations with $q$-expansions are limited by the computability of $q$-expansions of
modular forms, hence are generally not well suited for infinite families of computations.
In low degrees, however, our formula can directly be used to compute with
$q$-expansions. 
We demonstrate this by giving some sample calculations
of some $f$-invariants at the prime $5$.

\begin{remark}\label{rmk:23}
It is natural to ask if the results of this paper 
can be extended to the primes $2$ and
$3$.  A difficulty arises because the cohomology theory $\TMF$ fails to be
Landweber exact without inverting $6$, and this in turn is related to the
fact that the associated moduli stack of elliptic curves has geometric
points with automorphism groups divisible by the primes $2$ and $3$.  If
one
substitutes the group $SL_2(\ZZ)$ with a small enough congruence subgroup
so that the associated moduli stack is actually an algebraic space, then
the corresponding $f$-invariant detects the
$2$-line of the $2$ and $3$-primary Adams-Novikov spectral sequences.
However, the results of \cite{BM08} break down, because they rely on the
approximation theorem of \cite{BehrensLawson}, and the analog of this
approximation theorem for these congruence subgroups does not hold.  In fact,
the approximation theorem is not even true at the prime $2$ for the full
congruence subgroup $SL_2(\ZZ)$.
\end{remark}

We outline the organization of this paper.  In Section~\ref{sec:f}, we review the $f$-invariant.  
In Section~\ref{sec:f'}, we review the spectrum $Q(\ell)$, and use it to construct an
invariant $f'$ so that
$$ f_{i/j,k} = f'(\beta_{i/j,k}). $$
In Section~\ref{4} we show that the $f$-invariant is directly expressible in terms of
the invariant $f'$.  In Section~\ref{sec:examples}, we give our sample $5$-primary
calculations.

\section{The $f$-invariant}\label{sec:f}
This section reviews the $f$-invariant and its various aspects in
homotopy theory and geometry. Our main sources are
\cite{MR1781277} and \cite{MR1660325}.\par
\begin{thm}
Let $D$ be the ring
of divided congruences defined by N.\ Katz in \cite{MR0447119}, that
is, the ring of all inhomogeneous modular forms for $SL_2(\Z)$ whose
$q$-expansion is integral, and 
let $M_t$ be the subspace of modular forms of homogeneous weight $t$. 
Then for all $k>0$ there is a homomorphism
$$ f: \pi^ S_{2k}\lra D_\Q  /(D_{\Z[1/6]} \oplus (M_0)_\Q \oplus (M_{k+1})_\Q)
$$
whose kernel is the 3rd Adams-Novikov filtration for $MU[1/6]$.
\end{thm}

\begin{remark}
In \cite{MR1660325}, the second author actually defines the $f$ invariant
to take values in the subspace of 
$$ D_\Q  /(D_{\Z[1/6]} \oplus (M_0)_\Q \oplus (M_{k+1})_\Q) $$
spanned by inhomogeneous sums of modular forms of weights between $0$ and
$k+1$.  Of course, there is no harm in regarding the invariant as taking
values in
the larger group above.
\end{remark}

The construction of $f$ is closely related to the construction of the
classical $e$-invariant by F.\ Adams (see \cite{MR0198470}). Let $T$
be a flat ring spectrum  and let $$s:X\lra Y$$ be a stable map from a
finite spectrum into an arbitrary one. Suppose
further that
the $d$-invariant of $s$ vanishes. This simply means that $s$ vanishes
in $T$ homology. Then we have a short exact sequence
$$ {T}_* Y \lra {T}_ * C_s \lra {T}_* \Sigma X,$$
where $C_s$ is the cofiber of $s$. We can think of the sequence as an
extension of ${T}_*X$ by ${T}_*Y$ as a $T_*T$-comodule.
This is the classical $e$-invariant of $s$ in $T$-theory. \par
Next, suppose that
$$e(s)\in \Ext_{T_*T}({T}_*X,{T}_*Y)$$
vanishes, that is, the exact sequence of $T_*T$-comodules
splits and we choose a splitting. We also choose a
$T$-monomorphism $$ \iota: Y \lra I$$
into a $T$-injective spectrum $I$. For instance, we can take
$I=T\wedge Y$. Then there is a  map  $$t: C_s \lra I$$ which
is the image of $\iota_*$ under the induced splitting map
\[ [Y,I]\cong \Hom_{T_*T}({T}_*Y,{T}_*I)\lra \\ \Hom_{T_*T}
({T}_*C_s,{T}_*I)\cong [C_s,I].\]
In particular, the map $t$ coincides with $\iota$  on $Y$.
Let $F$ be the fiber of the map $\iota$. Then $s$ lifts to a map
$$  \bar{s}: X  \lra  F$$
which makes the diagram
$$ \xymatrix{\Sigma^{-1}C_s\ar[d]^{\Sigma^{-1} t} \ar[r]& X\ar[r]^s\ar[d]^{\bar{s}} & Y\ar[d]^{id}\\
\Sigma^{-1} I\ar[r] & F \ar[r]& Y }$$ commute.
\begin{lemma}\label{lem:dsbar}
$d(\bar{s})=0.$
\end{lemma}
\begin{proof}
In the split exact sequence
$$ \Hom_{T_*T}({T}_* \Sigma X,{T}_*\Sigma F)\lra
\Hom_{T_*T}({T}_* C_s,{T}_*\Sigma F)\lra
\Hom_{T_*T}({T}_* Y,{T}_*\Sigma F)
$$
the map $\Sigma \bar{s}_*$ restricted to $C_s$ is in the image of the
splitting and hence has to vanish. The claim follows since the map
from $C_s$ to $\Sigma X$ is surjective in $T$-homology.
\end{proof}
Lemma~\ref{lem:dsbar} implies that we again get a short exact sequence
$$ {T}_*F     \lra {T}_ * C_{\bar {s}} \lra {T}_ *
\Sigma X$$
which we can splice together with the short exact sequence
$$ {T}_ *\Sigma^ {-1} Y \lra  {T}_*\Sigma^ {-1} I
\lra {T}_ * F .$$
This gives an extension of ${T}_ *\Sigma^ {-1}Y$ by ${T}_ *\Sigma X$ of
length 2, that is, an element
$$ f(s)\in  \Ext^2_{T_ *T}({T}_ *X, {T}_* Y).$$
In the case $X=S^ {2k}$,  $Y=S^ 0$ and $T=\TMF[\tfrac{1}{6}]$, the image of
$f(s)$ under the 
injection
$$    \iota^2: \Ext^2 \hookrightarrow  
D_\Q  /(D_{\Z[\frac{1}{6}]} \oplus (M_0)_\Q \oplus (M_{k+1})_\Q) $$
is the second author's
$f$-invariant.  
The map $\iota^2$ will be reviewed in
Section~\ref{4}.
\par
We close this section with an alternative description of the
$f$-invariant. First recall from \cite{MR1781277} that a framed
manifold $M$ represents a framed bordism class  in second Adams-Novikov
filtration if and only if it is the corner of a $(U,fr)^ 2$ manifold
$W$. The boundary of $W$ is decomposed into two manifolds with
boundaries $W ^0$ and $W^ 1$. The stable tangent bundle of $W$ comes with a splitting
$$ TW \cong (TW)^{0} \oplus (TW)^{ 1}$$ and the bundles $(TW)^ i$ are
 trivialized on $W^i$.  Therefore, we get associated classes
$$ (TW)^ {i} \in K(W,W^ i).$$\par
Let $\exp_T$ be the usual
parameter for the universal Weierstrass cubic
$$ y^ 2 = 4 x^ 3 -E_4 x +E_6 $$
and let $$\exp_K(x)= 1-e^{- x}$$
be the standard parameter for the multiplicative formal group. Then
following theorem 
is a consequence of Proposition~4.1.4 of \cite{MR1781277} after
applying the complex orientation of the $\left < 2 \right>$-spectrum
$$ \xymatrix{S^ 0 \ar[r]\ar[d] & K \ar[d]\\
T \ar[r] & K \wedge T }.$$
\begin{thm} Let $s$ be represented by $M$ under the Pontryagin-Thom
  isomorphism. Then we have
$$f(s)= \left< \prod_{i,j} \frac{x_iy_j
  }{\exp_K(x_i)\exp_T(y_j)},[W,\partial W]\right>.$$
Here, $(x_i)$ and $(y_j)$ are the formal Chern roots of $(TW)^ 0$
and $(TW)^ 1$  respectively.
\end{thm}
We remark that there also are  descriptions of the $f$-invariant in
terms of a spectral invariant which is analogous to the classical
relation between the $e$-invariant and the $\eta$-invariant. We refer
the reader to \cite{hvb08} and \cite{BuNa08}.

\section{The spectrum $Q(\ell)$ and the invariant $f'$}\label{sec:f'}

For a $\ZZ[1/N]$-algebra $R$ we shall let
$M_k(\Gamma_0(N))_R$ denote the space of modular forms of weight $k$ over
$R$ of level $\Gamma_0(N)$ which are meromorphic at the cusps.  For $N = 1$
we shall simplify the notation by writing
$$ (M_k)_R := M_k(\Gamma_0(1))_R. $$

Let $\TMF_0(N)$ denote the corresponding spectrum of topological
modular forms with $N$ inverted (see \cite[Sec.~1.2.1]{BehrensK(2)},
\cite[Sec.~5]{Behrens}).  For primes $p > 3$, $\pi_* \TMF_0(N)_p$ is
concentrated in even degrees, and we have
\begin{equation}\label{eq:descent}
\pi_{2k} \TMF_0(N)_p \cong M_k(\Gamma_0(N))_{\ZZ_p}. 
\end{equation}

\begin{remark}
One could view the isomorphism of (\ref{eq:descent}) 
as a consequence of the fact that the
spectrum $\TMF_0(N)[\frac{1}{6}]$ is equivalent to the spectrum
$E^{\Gamma_0(N)}$ of \cite{MR1660325}, or as a consequence of the
fact that the descent spectral sequence
$$ H^s(\mc{M}^{\Gamma_0(N)}_{ell}[\tfrac{1}{6}], \omega^{\otimes t})
\Rightarrow \pi_{2t-s}
\TMF_0(N)[\tfrac{1}{6}] $$
is concentrated on $s=0$.
\end{remark}

Fix a pair of distinct primes $p$ and $\ell$.
In \cite{BehrensK(2)}, the first author introduced a $p$-local 
spectrum $Q(\ell)$,
defined as the totalization of a certain semi-cosimplicial spectrum
$$ Q(\ell) = \Tot(Q(\ell)^\bullet) $$
where $Q(\ell)^\bullet$ has the form
\begin{equation}\label{eq:cosimplicialQ}
Q(\ell)^\bullet =
\left( \TMF_p
\begin{array}{c}
\rightarrow \\
\rightarrow
\end{array}
\begin{array}{c}
\TMF_0(\ell)_p \\
\times \\
\TMF_p
\end{array}
\begin{array}{c}
\rightarrow \\
\rightarrow \\
\rightarrow
\end{array}
\TMF_0(\ell)_p
\right).
\end{equation}
In \cite[Sec.~4]{BM08} the spectrum
$Q(\ell)$ is reinterpreted as the smooth hypercohomology of a certain
open subgroup of an adele group 
acting on a certain spectrum.  
The semi-cosimplicial spectrum $Q(\ell)^\bullet$ is actually a
semi-cosimplicial $E_\infty$-ring spectrum, so the spectrum $Q(\ell)$ is an
$E_\infty$-ring spectrum.  In particular, there is a unit map
\begin{equation}\label{eq:unit}
\eta: S \rightarrow Q(\ell). 
\end{equation}
The spectrum $Q(\ell)$ is designed to be an
approximation of the $K(2)$-local sphere.  
More precisely, the spectrum
$Q(\ell)_{K(2)}$ is given as the homotopy fixed points of a subgroup
\begin{equation}\label{eq:Gammaell}
\Gamma_\ell \subset \MS_2
\end{equation}
of the Morava stabilizer group acting on the Morava $E$-theory $E_2$
\cite{Behrens} and this subgroup is dense if $\ell$ generates
$\ZZ_p^\times$ \cite{BehrensLawson}.   
The spectrum $Q(\ell)$ is
$E(2)$-local.
In \cite[Thm.~12.1]{BM08} it is proven that elements $\beta_{i/j,k} \in \pi_*(S_{E(2)})$ 
of \cite{MillerRavenelWilson} are detected by the map
$$ S_{E(2)} \rightarrow Q(\ell). $$
(It is not known if $Q(\ell)$ detects the entire divided beta family at the
primes $2$ and $3$.)

Taking the homotopy groups of the semi-cosimplicial spectrum $Q(\ell)^\bullet$
(\ref{eq:cosimplicialQ}) gives
a semi-cosimplicial abelian group
\begin{equation}\label{eq:cosimplicialgp}
C(\ell)_{2k}^\bullet := \left( (M_k)_{\ZZ_p}
\begin{array}{c}
\rightarrow \\
\rightarrow
\end{array}
\begin{array}{c}
M_k(\Gamma_0(\ell))_{\ZZ_p} \\
\times \\
(M_k)_{\ZZ_p}
\end{array}
\begin{array}{c}
\rightarrow \\
\rightarrow \\
\rightarrow
\end{array}
M_k(\Gamma_0(\ell))_{\ZZ_p}
\right).
\end{equation}
It is shown in \cite[Sec.~6]{BM08} that the 
morphisms
$$ d_0, d_1 : (M_k)_{\ZZ_p} \rightarrow M_k(\Gamma_0(\ell))_{\ZZ_p} \times
(M_k)_{\ZZ_p}, $$
induced by the initial coface maps of the cosimplicial abelian group
$C(\ell)^\bullet_{2k}$,
are given on the level of 
$q$-expansions by
\begin{align}
d_0(f(q)) & := (\ell^kf(q^\ell), \ell^kf(q)),  \label{eq:d0q} \\
d_1(f(q)) & := (f(q), f(q)). \label{eq:d1q}
\end{align}

The Bousfield-Kan spectral sequence for computing
$\pi_*\Tot(Q(\ell)^\bullet)$ gives a spectral sequence
\begin{equation}\label{eq:BKSS}
H^s(C(\ell)^\bullet)_t \Rightarrow \pi_{t-s} Q(\ell).
\end{equation}
For $p >3$, this spectral sequence collapses for dimensional reasons 
\cite[Cor.~5.2]{BM08}, giving us the following lemma.

\begin{lemma}\label{lem:BKSScollapse}
The edge homomorphism
$$ H^2(C(\ell)^\bullet)_t \rightarrow \pi_{t-2}(Q(\ell)) $$
is an isomorphism for $t \equiv 0 \mod 4$.
\end{lemma}

\begin{lemma}\label{lem:SSmap}
There is a map of spectral sequences
$$
\xymatrix{
\Ext^{s,t}_{BP_*BP}(BP_*, BP_*) \ar[d] \ar@{=>}[r] &
\pi_{t-s} S_{(p)} \ar[d]^{\eta_*} 
\\
H^s(C(\ell)^\bullet)_t \ar@{=>}[r] & 
\pi_{t-s}Q(\ell)
}
$$
from the Adams-Novikov spectral sequence for the sphere to the
Bousfield-Kan spectral sequence for $Q(\ell)$.
\end{lemma}

To prove Lemma~\ref{lem:SSmap} we shall
need the following lemma.

\begin{lemma}\label{lem:cosimpliciallemma}
Suppose that $R^\bullet$ is a semi-cosimplicial commutative $S$-algebra,
$E$ is a commutative $S$-algebra, and $\phi: E \rightarrow R^0$ is a map of 
commutative $S$-algebras.  Then there is a canonical extension of $\phi$ to a
map of semi-cosimplicial commutative $S$-algebras
$$ \phi^\bullet: E^{\wedge \bullet+1} \rightarrow R^\bullet $$
where
$$ E^{\wedge \bullet+1} = 
\left( E
\begin{array}{c}
\xrightarrow{\eta \wedge 1} \\
\xrightarrow{1 \wedge \eta} \\
\end{array}
E \wedge E
\begin{array}{c}
\xrightarrow{\eta \wedge 1 \wedge 1} \\
\xrightarrow{1 \wedge \eta \wedge 1} \\
\xrightarrow{1 \wedge 1 \wedge \eta}
\end{array}
E \wedge E \wedge E
\qquad \cdots 
\right)
$$
is the canonical cosimplicial $E$-resolution of the sphere.
\end{lemma}

\begin{proof}
A semi-cosimplicial commutative $S$-algebra is a functor
$$ \Delta_{\it inj} \rightarrow \{ \text{commutative $S$-algebras} \}, $$
where $\Delta_{\it inj}$ is the category of finite ordered sets and
order preserving injections.  Let $\ul{m}$ be the object of $\Delta_{\it
inj}$ given by 
$$ \ul{m} = \{ 0, 1, \ldots , m\} $$
and for $0 \le i \le n$ define $\iota^m_i: \ul{0} \rightarrow \ul{m}$ by
$\iota^m_i(0) = i$.  The map $\phi^s$ is defined to be the composite
$$ E^{\wedge s+1} \xrightarrow{((\iota_0^s)_* \circ \phi) \wedge \cdots
\wedge ((\iota_n^n)_* \circ \phi)} (R^s)^{\wedge s+1} \xrightarrow{\mu_{s+1}}
R^s $$
where $\mu_{s+1}$ denotes the $s+1$-fold product.  The maps $\phi^s$ are
easily seen to assemble into a map of semi-cosimplicial spectra.
\end{proof}

\begin{proof}[Proof of Lemma~\ref{lem:SSmap}]
Lemma~\ref{lem:cosimpliciallemma} implies that there exists a map of
semi-cosimplicial spectra
$$ 1^\bullet: \TMF_p^{\wedge \bullet + 1} \rightarrow Q(\ell)^\bullet $$
and hence a map from the Bousfield-Kan spectral sequence for
$\TMF_p^{\wedge \bullet + 1}$ to the Bousfield-Kan spectral sequence for
$Q(\ell)^\bullet$.  However, since $\TMF_p^{\bullet+1}$ is the canonical
$\TMF_p$-injective resolution of $S$, the Bousfield-Kan spectral sequence
for $\TMF_p^{\wedge \bullet + 1}$ is the $\TMF_p$-Adams-Novikov
spectral sequence for $S$.  Since $\TMF_p$ is complex orientable, there is
a map of ring spectra $BP \rightarrow \TMF_p$, and hence a map from the
$BP$-Adams-Novikov spectral sequence to the $\TMF_p$-Adams-Novikov spectral
sequence.
\end{proof}

The short exact sequences of $BP_*BP$-comodules
\begin{gather*}
0 \rightarrow BP_* \rightarrow BP_*[p^{-1}] \rightarrow BP_*/p^{\infty}
\rightarrow 0,
\\
0 \rightarrow BP_*/p^\infty \rightarrow BP_*/p^\infty[v_1^{-1}] \rightarrow 
BP_*/(p^\infty, v_1^\infty) \rightarrow 0
\end{gather*}
give rise to long exact sequences in $\Ext$, and the connecting homomorphisms
give a composite
\begin{multline}\label{eq:deltapv1}
\delta_{v_1,p}: \Ext^{0,t}_{BP_*BP}(BP_*, BP_*/(p^\infty, v_1^\infty)) \xrightarrow{\delta_{v_1}} 
\Ext_{BP_*BP}^{1,t}(BP_*, BP_*/p^\infty) 
\\
\xrightarrow{\delta_{p}}
\Ext_{BP_*BP}^{2,t}(BP_*, BP_*).
\end{multline}
The computations of \cite{MillerRavenelWilson} imply the following lemma.

\begin{lemma}\label{lem:deltapv1}
The homomorphism $\delta_{v_1,p}$ of (\ref{eq:deltapv1}) is an isomorphism
for $t > 0$.
\end{lemma}

Since the spectrum $\TMF[\frac{1}{6}]$ is Landweber exact, the spectrum
$\TMF_p$ is complex orientable.  Since $\TMF_p$ is $p$-local, it admits a
$p$-typical complex orientation, and 
a choice of $p$-typical complex orientation
$$ BP \rightarrow \TMF_p \rightarrow \TMF_0(\ell)_p $$
sends $v_1$ to a non-zero multiple of the Hasse invariant $E_{p-1}$ mod
$p$.  The complex $C(\ell)^\bullet/p^k$ is a complex of modules over the
ring $\ZZ_p[v_1^{p^{k-1}}]$. 
The short exact sequences
\begin{gather*}
0 \rightarrow C(\ell)^\bullet \rightarrow C(\ell)^\bullet[p^{-1}] \rightarrow 
C(\ell)^\bullet/p^{\infty}
\rightarrow 0,
\\
0 \rightarrow C(\ell)^\bullet/p^\infty \rightarrow 
C(\ell)^\bullet/p^\infty[v_1^{-1}] \rightarrow 
C(\ell)^\bullet/(p^\infty, v_1^\infty) \rightarrow 0
\end{gather*}
give rise to long exact sequences in $H^*$, and the connecting homomorphisms
give a composite
\begin{equation*}
\delta_{v_1,p}: H^{0}(C(\ell)^\bullet/(p^\infty, v_1^\infty))_t 
\xrightarrow{\delta_{v_1}} 
H^{1}(C(\ell)^\bullet/p^\infty)_t \xrightarrow{\delta_{p}}
H^{2}(C(\ell)^\bullet)_t.
\end{equation*}

Using Lemmas~\ref{lem:BKSScollapse} and \ref{lem:deltapv1}, we have the
following diagram, for $t > 0$.
\begin{equation}\label{diag:f'}
\xymatrix{
\pi_{4t-2} S_{(p)} \ar[r] &
\pi_{4t-2} Q(\ell)
\\
\Ext^{2,4t}_{BP_* BP_*}(BP_*, BP_*) \ar[r] \ar@{=>}[u] \ar@{.>}[dr]^{f'} &
H^{2}(C(\ell)^\bullet)_{4t} \ar[u]_\cong
\\
\Ext^{0,4t}_{BP_*BP}(BP_*, BP_*/(p^\infty, v_1^\infty))
\ar[u]^{\cong}_{\delta_{v_1,p}} \ar[r] &
H^{0}(C(\ell)^\bullet/(p^\infty, v_1^\infty))_{4t} \ar[u]_{\delta_{v_1,p}}
}
\end{equation}
Since $p$ is odd and $\Ext^{2,m}_{BP_* BP}(BP_*, BP_*)$ is concentrated
in degrees $m \equiv 0 \mod 4$, the invariant $f'$ may be regarded as an
invariant defined on the entire $2$-line of the ANSS.
Moreover, because $\pi_{4t-2} S_{(p)}$ 
contains no elements of Adams-Novikov filtration less
than $2$, the invariant $f'$ may be regarded as giving a homotopy invariant
through the composite
$$ \pi_{4t-2} S_{(p)} \rightarrow \Ext^{2,4t}_{BP_*BP}(BP_*, BP_*) 
\xrightarrow{f'}
H^{0}(C(\ell)^\bullet/(p^\infty, v_1^\infty))_{4t}. $$
We shall find that this invariant $f'$ is closely related to the $f$
invariant of the second author.

We end this section by describing some of the salient features of the 
invariant $f'$.  Namely, we shall show:
\begin{enumerate}
\item  
the homomorphism $f'$ is a monomorphism, and if $\ell$ generates
$\ZZ_p^\times$, the homomorphism $f'$ is almost an
isomorphism, and 
\item the groups $H^0(C(\ell)^\bullet/(p^\infty, v_1^\infty))_{4t}$ 
admit a precise arithmetic
interpretation in terms of congruences of $q$-expansions of modular forms. 
\end{enumerate}

\subsection*{The injectivity and almost surjectivity of $f'$}
$\quad$

Because $v_2$ is invertible in $C(\ell)^\bullet/(p^\infty, v_1^\infty)$,
there is a factorization
\begin{equation}\label{diag:f'factor}
\xymatrix{
\Ext^{2,4t}_{BP_*BP}(BP_*, BP_*) \ar[r]^{f'} &
H^0(C(\ell)^\bullet/(p^\infty, v_1^\infty))_{4t} 
\\
\Ext^{0,4t}_{BP_*BP}(BP_*, BP_*/(p^\infty, v_1^\infty)) \ar[u]^{\cong}_{\delta_{v_1,p}} \ar[r]_{L_{v_2}} &
\Ext^{0,4t}_{BP_*BP}(BP_*, BP_*/(p^\infty, v_1^\infty)[v_2^{-1}])
\ar@{.>}[u]_{\bar{\eta}}
}
\end{equation}

Recall from \cite{MillerRavenelWilson} that for $t > 0$ the groups
$$ \Ext^{0,4t}_{BP_*BP}(BP_*, BP_*/(p^\infty, v_1^\infty)) 
\quad \text{and} \quad 
\Ext^{0,4t}_{BP_*BP}(BP_*, BP_*/(p^\infty, v_1^\infty)[v_2^{-1}])
$$
are generated by elements $\beta_{i/j,k}$ for certain combinations of
indices $i$, $j$, and $k$.  As usual, $\beta_{i/j}$ denotes the element
$\beta_{i/j,1}$.

\begin{prop}\label{prop:etabar}
$\quad$
\begin{enumerate}
\item The map $L_{v_2}$ of (\ref{diag:f'factor}) is injective, and the cokernel is an $\FF_p$-vector
space with basis
$$ \{ \beta_{p^n/j} \: : \: n \ge 2, p^n < j \le p^n + p^{n-1} - 1 \}. $$

\item The map $\bar{\eta}$ of (\ref{diag:f'factor}) is injective, and if
$\ell$ generates $\ZZ_p^\times$, it is an isomorphism.
\end{enumerate}
\end{prop}

\begin{proof}
(i) follows directly from the calculations of \cite{MillerRavenelWilson}.
(ii) follows from the fact that the map $\bar{\eta}$ factors as
$$
\xymatrix@C-5em@R+1em{
\Ext_{BP_*BP}^{0,4t}(BP_*, BP_*/(p^\infty, v_1^\infty)[v_2^{-1}])
\ar[rr]^-{\bar{\eta}} \ar[dr]^{\bar{\eta}'}_{\cong} && H^0(C(\ell)^\bullet/(p^\infty,
v_1^\infty))_{4t} 
\\
& H^0_c(\MS_2, \pi_{*}E_2/(p^\infty, v_1^\infty))^{\Gal}_{4t}
\ar[ur]^{\bar{\eta}''}
}
$$
where $\bar{\eta}'$ is the Morava change-of-rings isomorphism, and
$\bar{\eta}''$ is given by the composite
$$ H^0(C(\ell)^\bullet/(p^\infty, v_1^\infty))_{4t}
\xrightarrow[\cong]{\omega}  
H^0(\Gamma_\ell, \pi_{*}M_2 E_2)^{\Gal}_{4t}
\xrightarrow{\nu}
H^0_c(\MS_2, \pi_{*}M_2 E_2)^{\Gal}_{4t}.
$$
Here, $\omega$ is the isomorphism given by
\cite[Cor.~7.7]{BM08} where $\Gamma_\ell$ is the subgroup of $\MS_2$ of
(\ref{eq:Gammaell}), the
spectrum $M_2 E_2$ is the second monochromatic layer of $E_2$, and
$\nu$ is the monomorphism induced by the inclusion of the subgroup.
Lemma~11.1 of \cite{BM08} states that $\nu$ is an 
isomorphism if $\ell$ generates
$\ZZ_p^\times$ \cite[Lem.~11.1]{BM08}.  Note that the same argument in
\cite[Thm.~6.1]{Ravenel84} computing $\pi_* M_n BP$ applies to compute
$$ \pi_* M_2 E_2 \cong (\pi_{*}E_2)/(p^\infty, v_1^\infty). $$
\end{proof}

We conclude that $f'$ is injective, and if $\ell$ generates $\ZZ_p^\times$, 
the only generators of
$H^0(C(\ell)^\bullet/(p^\infty, v_1^\infty))$
not in the
image of $f'$ are those corresponding to the Greek letter elements 
$\beta_{p^n/j}$ for $j > p^n$.

\subsection*{The arithmetic interpretation of the groups
$H^0(C(\ell)^\bullet/(p^\infty, v_1^\infty))$} $\quad$

The groups $H^0(C(\ell)^\bullet/(p^\infty, v_1^\infty))_{4t}$ 
are computed by the colimit of
groups
$$ H^0(C(\ell)^\bullet/(p^\infty, v_1^\infty))_{4t} = \colim_{k} \colim_{\substack{j = sp^{k-1} \\
s \ge 1}} \mc{B}_{2t/j,k} $$
where
$$ \mc{B}_{2t/j,k} = H^0(C(\ell)^\bullet/(p^k, v_1^j))_{4t + 2j(p-1)}. $$
Using the fact that $v_1$ corresponds, modulo $p$, 
to a non-zero multiple of the 
Hasse invariant $E_{p-1}$ in
the ring of modular forms, we have
$$
\mc{B}_{2t/j,k}
= \ker \left(
\frac{M_{2t+j(p-1)}}{(p^k, E_{p-1}^j)} \xrightarrow{d_0 - d_1}
\begin{array}{c}
\frac{M_{2t+j(p-1)}}{(p^k, E_{p-1}^j)} \\
\oplus \\
\frac{M_{2t+j(p-1)}(\Gamma_0(\ell))}{(p^k, E_{p-1}^j)}
\end{array}
\right).
$$

Serre \cite[Prop.~4.4.2]{MR0447119} 
showed that two modular forms $f_1$ and $f_2$ 
over $\ZZ/p^k$ are linked by
multiplication by $E_{p-1}^j$
(for $j \equiv 0 \mod p^{k-1}$) if and only if the corresponding
$q$-expansions satisfy 
$$ f_1(q) \equiv f_2(q) \mod p^k. $$ 
Using this, and (\ref{eq:d0q})-(\ref{eq:d1q}), the following theorem is
proven in \cite{BM08}.

\begin{thm}[{\cite[Thm.~11.3]{BM08}}]\label{thm:BMthm}
There is a one-to-one correspondence between the additive generators of
order $p^k$ in $\mc{B}_{t/j,k}$
and the modular forms $f \in M_{t+j(p-1)}$ (modulo $p^k$) satisfying
\begin{enumerate}
\item[$\rm (1)$] We have $t \equiv 0 \mod (p-1)p^{k-1}$.
\item[$\rm (2)$] The $q$-expansion $f(q)$ is not congruent to $0$ mod $p$.
\item[$\rm (3)$] We have $\ord_q f(q) > \frac{t}{12}$ or $\ord_q f(q) =
\frac{t-2}{12}$.
\item[$\rm (4)$] There does not exist a form $f' \in M_{t'}$ such that
$f'(q) \equiv f(q) \mod p^k$ for $t' < t + j(p-1)$.
\item[$\rm{(5)}_\ell$] There exists a form 
$$ g \in M_{t}(\Gamma_0(\ell)) $$
satisfying 
$$ f(q^\ell) - f(q) \equiv  g(q) \mod p^k. $$
\end{enumerate}
\end{thm}

\begin{remark}
It follows from \cite[Cor.~11.7]{BM08}, that 
a modular form satisfying (1)--(5) corresponding to $f'(x)$ is independent
of the choice of the prime $\ell$. 
\end{remark}

\section{The relation between $f$ and $f'$}\label{4}
Let $\ell$ be a generator of $\ZZ_p^\times$. 
We start with a cohomology class
$$ x \in \Ext^{2,2t}_{BP_*BP}(BP_*, BP_*) $$
with corresponding invariant
\begin{equation}\label{eq:f'jk}
f'(x) \in \mathcal{B}_{t/j,k}=
H^ 0(C^ \bullet(\ell)/(p^k,v_1^j))_{2t+2j(p-1)}.
\end{equation}
Note that since $p$ is odd, $t$ must be even.
By Theorem~\ref{thm:BMthm}, a 
representative of $f'(x)$ is a $\Z/p^ k$ modular form $\varphi$ 
of weight $t+j(p-1)$
for $SL_2(\Z)$ which satisfies
certain congruences.  We view $\varphi$ as a divided
congruence, more precisely, as an element of
\[ D\otimes \Z /p^ k .\]

\begin{thm}\label{thm:mainthm}
The $f$-invariant of the class $x$ is given by
\[ p^{-k}E_{p-1}^{-j}(\varphi-q^ 0(\varphi))\]
where $q^ 0$ is the 0th Fourier coefficient, and $j,k$ are given by
(\ref{eq:f'jk}).
\end{thm}

The proof of Theorem~\ref{thm:mainthm} will be deferred to the end of the
section.

\begin{remark}
For $t > 0$, Theorem~\ref{thm:BMthm}(3) implies that there exists a
representative $\varphi$ of $f'(x)$ with $q^0(\varphi) = 0$.  Since the
modular form $f_{i/j,k}$ of \cite{BM08} is such a representative of 
$f'(\beta_{i/j,k})$, Theorem~\ref{thm:mainthm} implies that
$$ f(\beta_{i/j,k}) = \frac{f_{i/j,k}}{p^k E_{p-1}^j}. $$
\end{remark}

\begin{cor}
The class
\[  p^ k E_{p-1}^j f(x)\]
is congruent to  a $\Z/p^k$-modular form $\varphi$ of weight $t+j(p-1)$
up to modular forms of weights $j(p-1)$ and $t+j(p-1)$. Moreover,
$\varphi$ satisfies the conditions (1)-(5) of \ref{thm:BMthm}.
\end{cor}

\begin{remark}
We pause to explain how the expression in Theorem~\ref{thm:mainthm} may be
regarded as an element of the subgroup
$$ \frac{D_\QQ}{D_{\ZZ_{(p)}} + (M_0)_\QQ + (M_t)_\QQ}
\subset
\frac{D_\QQ}{D_{\ZZ[1/6]} + (M_0)_\QQ + (M_t)_\QQ} $$
in a way that more clearly accounts for the indeterminacy of the $f$-invariant.
Katz showed that $D$ is a dense subspace of $\mb{V}$, the ring of
generalized $p$-adic modular functions \cite{Katz}.  The ring $\mb{V}$ has
an action by the group $\ZZ_p^\times$ through Diamond operators, and the
weight $t$ subspace $\mb{V}_t$ is canonically identified by
$$ \mb{V}_t \cong (M_*)_{\ZZ_p}[E_{p-1}^{-1}]_t. $$  
We therefore have
\begin{align*}
\frac{D_\QQ}{D_{\ZZ_{(p)}} + (M_0)_\QQ + (M_t)_\QQ}
& \cong \frac{\mb{V}_\QQ}{\mb{V} + (M_0)_{\QQ_p} + (M_t)_{\QQ_p}}.
\end{align*}
Taking the weight $t$ subspace we get
\begin{align*} 
\frac{(\mb{V}_t)_\QQ}{\mb{V}_t + (M_t)_{\QQ_p}} 
& \cong \left( \frac{(M_*)_{\QQ_p}[E_{p-1}^{-1}]}{(M_*)_{\ZZ_p}[E_{p-1}^{-1}] 
+ (M_*)_{\QQ_p}} \right)_t \\
& = \left( \frac{(M_*)_{\ZZ_p}}{(p^\infty, E_{p-1}^\infty)} \right)_t.
\end{align*}
The expression $p^{-k}E_{p-1}^{-j}\phi$ clearly may be regarded as an 
element of the group above.
\end{remark}

Let $T$ be $\TMF[\tfrac{1}{6}]$ and 
$$ M^ {(2)} = \pi_* T \wedge T $$ 
be the Hopf algebroid of 
cooperations
of $T$. An element of $M^ {(2)}$ is a modular form in two variables
which is meromorphic at $\infty$ and has (away from 6) an integral
Fourier expansion (see \cite{MR1660325}). 

Consider the map of semi-cosimplicial spectra
$$ 1^\bullet: \TMF_p^{\wedge \bullet + 1} \rightarrow Q(\ell)^\bullet $$
of
Lemma~\ref{lem:cosimpliciallemma}. 
Applying the functor $\pi_*(-)$, we get a map of semi-cosimplicial abelian
groups 
$$ \pi_{2k}(T_p^{\bullet+1}) = M_k^{(\bullet + 1)} \rightarrow
C^\bullet(\ell)_{2k} $$
which in low
degrees gives the following commutative diagram.
$$ \xymatrix{(M_k)_{\Z_p}\ar[rr]^ {d_0-d_1} \ar[d]^ = &&
(M^ {(2)}_k)_{\Z_p}\ar[d]^ \phi\\
(M_k)_{\Z_p}\ar[rr]^ -{d_0-d_1}&&
  M_k(\Gamma_0(\ell))_{\Z_p}\times
(M_k)_{\Z_p}}.$$

\begin{lemma}
The induced map in cohomology
$$ H^ 0(M^{(\bullet+1)}_*/(p^ \infty,E_{p-1}^ \infty))\lra 
H^ 0(C^ \bullet(\ell)/(p^ \infty,v_1^\infty))$$
is an isomorphism.
\end{lemma}

\begin{proof}
By \cite{HoveyStrickland}, there is a change-of-rings isomorphism
\begin{align*}
H^ 0(M^{(\bullet+1)}_*/(p^ \infty,E_{p-1}^ \infty)) 
& = \Ext^0_{\TMF_*\TMF_p}(\pi_* \TMF_p,
\pi_* \TMF_p/(p^\infty, E_{p-1}^\infty)) \\
& \cong \Ext^0_{BP_*BP}(BP_*, BP_*/(p^\infty, v_1^\infty)[v_2^{-1}]).
\end{align*}
The lemma follows from the isomorphism $\bar{\eta}$ of
Proposition~\ref{prop:etabar}.
\end{proof}

Next we explain how to get from an element in
$$  H^ 0(M_*/(p^ \infty,E_{p-1}^ \infty))\cong \Ext^ 0_{M^ {(2)}}(M_*,M_*/(p^ \infty,E_{p-1}^ \infty))$$
to a congruence in
$$  D_\Q  /(D_{\Z[1/6]} \oplus (M_0)_\Q \oplus (M_{k})_\Q).$$
For this, we first describe how a class $\varphi$ in
$$ \Ext^ 0_{M^ {(2)}}(M_*,M_*/(p^ \infty,E_{p-1}^ \infty))$$
gives rise to a class in
$$ \Ext^ 2 _{M^ {(2)}}(M_*,M_*).$$
We use the geometric boundary theorem
\begin{thm}\label{GBT}
\cite{MR860042}
Write $E_*(X)$ for the $E_*$-term of the $T$-based Adams Novikov
spectral sequence which conditionally converges to the homotopy of the
$T$-nilpotent
completion of $X$.  Let
\[
W\stackrel{f}{\lra}X\stackrel{g}{\lra}Y \stackrel{h}{\lra}\Sigma W
\]
be a cofiber sequence of finite spectra with $T_*( h)=0$.
Assume further that $[s]\in {E}_2^{t,*+t}(Y)$ converges to
$s$. Then $\delta \,  {[s]}$ converges to
$h_*(s)$ where $\delta$ is the connecting
homomorphism to the short exact sequence of chain complexes
\[
0\lra {E}_1(W)\lra {E}_1(X)\lra {E}_1(Y)\lra 0.
\]
\end{thm}
For a multi index $I$ let $$M(I)=M(i_0,\ldots ,i_{n-1})$$ be the
generalized Moore spectrum with
$$BP_*M(I)= \Sigma ^{-||I||-n}BP_*/(p^ {i_0},v_1^ {i_1}, \ldots,
  v_{n-1}^ {i_{n-1}})$$
where
$$ ||I||=\sum_j2i_j(p^ j-1)$$
Each $M(I)$ admits a self map
$$ \Sigma^ {2i_n(p^ n-1)} M(I) \lra M(I)$$
which induces multiplication by $v_n^ {i_n}$. Its fiber is
$M(I,i_n)$. We apply the geometric boundary theorem to the sequences
$$ \xymatrix{ \Sigma^ {2i_1(p-1)} M(i_0) \ar[r]^-{v_1^i} & M(i_0)
\ar[r] &
\Sigma M(i_0,i_1)\ar[r] &\Sigma^ {2i_1(p-1)+1} M(i_0)}  $$
and
$$ \xymatrix{ S \ar[r]^{p^ {i_0}}& S\ar[r] &  \Sigma M(i_0)\ar[r]& S^ 1
    }.$$
For $$
\varphi \in E^ {0}_2(M(i_0,i_1))=\Ext^ 0_{M^ {(2)}}(M_*,M_*/(p^
  {i_0},E_{p-1}^ {i_1}))$$
we have
$$ \delta \varphi = \left[ \frac{d^ 0\varphi-d^ 1 \varphi}{E_{p-1}^{i_1}}
\right] \in
E^ {1}_2(M(i_0))=\Ext^ 1_{M^ {(2)}}(M_*,M_*/p^{i_0})$$
and
$$ \delta \delta \varphi =\left[ p^ {-i_0}\sum_{i=0}^ 2 (-1)^
  i d^ {i} \left[ \frac{d^ 0\varphi-d^ 1 \varphi}{E_{p-1}^{i_1}}
  \right]\right]\in
E^ {2}_2(S)=\Ext^ 2_{M^ {(2)}}(M_*,M_*) $$
where $d^ i$ denote the differentials of the cobar complex
$$ (\Omega_T^\bullet)_{2k} = \pi_{2k} T^{\bullet+1} \cong M_k^{(\bullet+1)}. $$ 
The maps of ring spectra
$$ \xymatrix{T \ar[r]^ {q^ 0} & K_{\Z[1/6]} \ar[r]^ {ch^ 0}& H_{\Q}}$$
induce the following map of semi-cosimplicial spectra.
$$ 
\xymatrix@C+6em@R+2em{
T
\ar@<1ex>[r]|{\eta \wedge 1}
\ar@<-1ex>[r]|{1 \wedge \eta}
\ar[d]^{1}
&
T \wedge T
\ar@<2ex>[r]|{\eta \wedge 1 \wedge 1}
\ar[r]|{1 \wedge \eta \wedge 1}
\ar@<-2ex>[r]|{1 \wedge 1 \wedge \eta}
\ar[d]|{q^0 \wedge 1}
&
T \wedge T \wedge T
\ar[d]|{\mathit{ch}^0 \circ q^0 \wedge q^0 \wedge 1}
\\
T
\ar@<1ex>[r]|{\eta \wedge 1}
\ar@<-1ex>[r]|{q^0 \wedge \eta}
&
K_{\ZZ[1/6]} \wedge T
\ar@<2ex>[r]|{\eta \wedge 1 \wedge 1}
\ar[r]|{\mathit{ch}^0 \wedge \eta \wedge 1}
\ar@<-2ex>[r]|{\mathit{ch}^0 \wedge q^0 \wedge \eta}
&
H_\QQ \wedge K_{\ZZ[1/6]} \wedge T
} $$
Taking $\pi_{2k}(-)$, and using 
\cite[Thm.~2.7]{MR1660325}, we get the following map of semi-cosimplicial abelian groups
$$ 
\xymatrix@C+0em@R+1em{
(\Omega_T^\bullet)_{2k} \ar[d]_\rho
&
M_k^{(1)}
\ar@<1ex>[r]|{d_0}
\ar@<-1ex>[r]|{d_1}
\ar[d]_=^{\rho^0}
&
M^{(2)}_k
\ar@<2ex>[r]|{d_0}
\ar[r]|{d_1}
\ar@<-2ex>[r]|{d_2}
\ar[d]^{\rho^1}
&
M^{(3)}_k
\ar[d]^{\rho^2}
\\
(\Omega^\bullet_{T,K,H})_{2k}
&
(M_k)_{\ZZ[\tfrac{1}{6}]}
\ar@<1ex>[r]|{d_0}
\ar@<-1ex>[r]|{d_1}
&
D_{\ZZ[\tfrac{1}{6}]}
\ar@<2ex>[r]|{d_0}
\ar[r]|{d_1}
\ar@<-2ex>[r]|{d_2}
&
D_\QQ
} $$
In $(\Omega^\bullet_{T,K,H})_{2k}$, we have
\begin{gather*}
d_1(D_{\ZZ[\tfrac{1}{6}]}) \subseteq (M_k)_\QQ \subseteq D_\QQ, \\
d_2(D_{\ZZ[\tfrac{1}{6}]}) \subseteq (M_0)_\QQ \subseteq D_\QQ.
\end{gather*}
Therefore, by modding out by these subgroups of $D_\QQ$, we get a map:
$$ 
\xymatrix@C+0em@R+1em{
(\Omega_T^\bullet)_{2k} \ar[d]_{\bar\rho}
&
M_k^{(1)}
\ar@<1ex>[r]|{d_0}
\ar@<-1ex>[r]|{d_1}
\ar[d]_=^{\bar\rho^0}
&
M^{(2)}_k
\ar@<2ex>[r]|{d_0}
\ar[r]|{d_1}
\ar@<-2ex>[r]|{d_2}
\ar[d]^{\bar\rho^1}
&
M^{(3)}_k
\ar[d]^{\bar\rho^2}
\\
(\bar{\Omega}^\bullet_{T,K,H})_{2k}
&
(M_k)_{\ZZ[\tfrac{1}{6}]}
\ar@<1ex>[r]|{d_0}
\ar@<-1ex>[r]|{d_1}
&
D_{\ZZ[\tfrac{1}{6}]}
\ar@<2ex>[r]|-{d_0}
\ar[r]|-{d_1}
\ar@<-2ex>[r]|-{d_2}
&
D_\QQ/((M_0)_\QQ + (M_k)_\QQ)
} $$
The
first coface maps of the semi-cosimplicial abelian group
$(\bar{\Omega}^\bullet_{T,K,H})_{2k}$ are given by
$$ d^ 0 = \iota,\;  d^ 1 =q^ 0$$
and the second ones by
$$d^ 0= \iota, \; d^ 1 = d^ 2 =0 $$
where $\iota$ the canonical inclusion. 
The induced map in cohomology
is the inclusion
$$\iota^2 : \Ext^ 2 _{M^ {(2)}}(M_*,M_*) \hookrightarrow  D_\Q  /((D)_{\Z[1/6]}
  \oplus (M_0)_\Q \oplus (M_{k})_\Q).$$
Hence we have
$$ \bar\rho_*\delta \delta \varphi =  p^{-i_0}E_{p-1}^{-i_1}(\varphi-q^
0(\varphi))$$
and the proof of the theorem is completed.

\section{Examples at $p = 5$}\label{sec:examples}

Below are some computations of the $q$-expansions of the modular forms 
$f_{i/j,k}$
representing $f'(\beta_{i/j,k})$ at $p = 5$.  
The $q$-expansions of the 
corresponding $f$ invariants, by
Theorem~\ref{thm:mainthm}, are given by 
$$ f(\beta_{i/j,k}) = p^{-k} E^ {-j}_{p-1} f_{i/j,k}(q). $$
The computations were performed using the MAGMA computer algebra system, with $\ell = 2$,  
as follows.
\begin{enumerate}
\item A basis $\{F_\alpha(q) \}$ of $q$-expansions of forms in $M_{24i}$ satisfying
Theorem~\ref{thm:BMthm}(3) was generated.
\item A basis $\{ G_\beta(q) \}$ of $q$-expansions of holomorphic forms in 
$M_{24i-4j}(\Gamma_0(\ell))_{\ZZ/5^k}$ was generated.
\item Basic linear algebra is used to calculate a basis of linear combinations
$ \sum_\alpha a_\alpha F_\alpha $
such that 
$$ \sum_\alpha a_\alpha (F_\alpha(q^2) - F_\alpha(q)) \equiv 
\sum_\beta b_\beta G_\beta(q) \mod 5^k. $$
\end{enumerate}

\begin{note}
The following modular forms are normalized so that the leading term has coefficient $1$.
Therefore, they may differ from the $f'$-invariants of $\beta_{i/j,k}$ by a multiple in
$\ZZ_p^\times$.
\end{note}

\begin{supertabular}{rp{4in}}
$f_{1/1,1} =$ &
$\Delta^2 =$ 
\\
&
\begin{verbatim}q^2 + 2*q^3 + q^7 + q^12 + 2*q^13 + q^17 + 2*q^18 + 
2*q^22 + 2*q^23 + 3*q^28 + q^32 + 4*q^33 + q^37 + 
2*q^42 + 2*q^43 + q^47 + 2*q^48 + q^52 + 2*q^53 + 
2*q^62 + 2*q^63 + q^67 + 3*q^68 + 2*q^73 + 2*q^77 + 
4*q^78 + 2*q^82 + 2*q^83 + q^92 + 4*q^93 + q^97 + 
3*q^98 + O(q^100) mod 5\end{verbatim}
\\
$f_{2/1,1} =$ & $\Delta^4 =$ 
\\
&
\begin{verbatim}q^4 + 4*q^5 + 4*q^6 + 2*q^9 + 4*q^10 + 3*q^14 + 
3*q^15 + 3*q^16 + 4*q^19 + 2*q^20 + 3*q^21 + 2*q^24 
+ 2*q^26 + q^29 + 3*q^30 + 2*q^34 + 4*q^35 + 3*q^36
+ 3*q^39 + 2*q^44 + 3*q^45 + q^51 + 4*q^54 + 3*q^55 
+ q^56 + 2*q^59 + 4*q^60 + 2*q^64 + 3*q^65 + 3*q^66 
+ 4*q^69 + 4*q^70 + 2*q^76 + q^79 + 4*q^80 + 4*q^81 
+ q^84 + 4*q^85 + q^86 + 3*q^89 + 3*q^90 + q^91 + 
4*q^94 + 4*q^96 + 4*q^99 + O(q^100) mod 5\end{verbatim}
\\
$f_{3/1,1} =$ & $\Delta^6 = $ 
\\
&
\begin{verbatim}q^6 + q^7 + 2*q^8 + 3*q^9 + 3*q^11 + 2*q^12 + 2*q^13 
+ q^16 + 4*q^17 + q^18 + 4*q^19 + 2*q^22 + 4*q^24 +
3*q^26 + 3*q^27 + 3*q^28 + 3*q^29 + 4*q^31 + 4*q^32 
+ 4*q^33 + 4*q^34 + q^36 + q^37 + 4*q^38 + 3*q^39 + 
4*q^41 + q^42 + 4*q^44 + 4*q^46 + 4*q^48 + 4*q^49 + 
q^51 + 2*q^53 + 4*q^54 + 3*q^56 + 4*q^58 + q^62 + 
4*q^63 + 3*q^64 + 3*q^66 + 4*q^67 + 3*q^68 + q^69 + 
2*q^72 + 4*q^73 + q^74 + q^76 + 4*q^77 + 3*q^78 + 
4*q^79 + q^82 + 3*q^84 + 2*q^86 + q^87 + 4*q^88 + 
4*q^89 + 3*q^91 + q^92 + 2*q^93 + 4*q^94 + 3*q^96 + 
3*q^97 + q^98 + 2*q^99 + O(q^100) mod 5\end{verbatim}
\\
$f_{4/1,1} =$ & $\Delta^8 = $ 
\\
&
\begin{verbatim}q^8 + 3*q^9 + 4*q^10 + 2*q^11 + q^12 + 4*q^13 + 4*q^14
+ 3*q^15 + 2*q^16 + q^19 + 3*q^21 + 4*q^22 + 2*q^24 
+ 4*q^26 + 4*q^27 + 4*q^28 + 4*q^29 + 3*q^31 + 
4*q^33 + q^34 + 4*q^35 + 3*q^37 + q^38 + 2*q^39 + 
q^43 + 3*q^44 + 2*q^47 + 4*q^51 + 2*q^52 + q^53 + 
3*q^54 + q^56 + q^57 + 3*q^58 + 2*q^59 + 4*q^60 +
4*q^61 + 2*q^63 + 3*q^65 + 2*q^66 + q^67 + 4*q^68 + 
2*q^69 + 2*q^71 + q^73 + q^74 + 2*q^76 + 2*q^78 + 
3*q^79 + 2*q^81 + 3*q^82 + 4*q^85 + 4*q^86 + q^87 +
q^89 + 3*q^90 + q^91 + 3*q^92 + 3*q^93 + 3*q^94 + 
4*q^97 + 3*q^98 + 4*q^99 + O(q^100) mod 5\end{verbatim}
\\
$f_{5/5,1} =$ &
$\Delta^{10} =$ 
\\
&
\begin{verbatim}q^10 + 2*q^15 + q^35 + q^60 + 2*q^65 + q^85 + 2*q^90 
+ O(q^100) mod 5\end{verbatim}
\\
$f_{25/29,1} =$ & 
$\Delta^{50} + 4\Delta^{42}E_4^{24} + 3\Delta^{41}E_4^{27} = $
\\
&
\begin{verbatim}3*q^41 + 2*q^42 + 4*q^43 + 4*q^44 + 3*q^47 + 2*q^48 + 
3*q^49 + q^50 + q^51 + q^52 + 2*q^54 + q^56 + 4*q^58
+ q^59 + 4*q^61 + 4*q^62 + q^63 + 3*q^64 + q^66 + 
4*q^67 + 3*q^68 + 3*q^69 + q^71 + q^74 + 2*q^75 + 
2*q^76 + 3*q^78 + 4*q^79 + 2*q^81 + 3*q^82 + 2*q^83 
+ 4*q^84 + 2*q^88 + 3*q^89 + 4*q^91 + q^92 + 2*q^94 
+ 2*q^96 + q^98 + q^102 + q^104 + 4*q^106 + 3*q^107 
+ 3*q^108 + 2*q^109 + 4*q^111 + 4*q^112 + 4*q^114 + 
3*q^116 + 2*q^118 + 2*q^119 + q^121 + 4*q^122 + 
3*q^123 + q^124 + q^126 + 2*q^127 + q^129 + 4*q^132 
+ q^134 + 4*q^136 + 4*q^138 + q^139 + q^141 + 
3*q^143 + q^144 + q^147 + 3*q^149 + O(q^150) mod 5\end{verbatim}
\\
$f_{25/5,2} =$ &
$\Delta^{50} = $ \\ &
\begin{verbatim}q^50 + 10*q^55 + 15*q^60 + 5*q^65 + 5*q^70 + 12*q^75 + 
15*q^80 + 20*q^85 + 10*q^90 + 5*q^95 + 15*q^100 + 
10*q^105 + 20*q^110 + 5*q^115 + 20*q^125 + 20*q^135 
+ 15*q^140 + 20*q^145 + 10*q^150 + O(q^151) mod 25\end{verbatim}
\\
\end{supertabular}

\bibliographystyle{amsalpha}
\bibliography{qf3}
\end{document}